\input amstex
 \documentstyle{amsppt}
 \magnification=\magstep1
 \vsize=24.2true cm
 \hsize=15.3true cm
 \nopagenumbers\topskip=1truecm
 \headline={\tenrm\hfil\folio\hfil}

 \TagsOnRight

\hyphenation{auto-mor-phism auto-mor-phisms co-homo-log-i-cal
co-homo-logy
co-homo-logous dual-izing pre-dual-izing geo-metric geo-metries
geo-metry
half-space homeo-mor-phic homeo-mor-phism homeo-mor-phisms
homo-log-i-cal
homo-logy homo-logous homo-mor-phism homo-mor-phisms hyper-plane
hyper-planes
hyper-sur-face hyper-sur-faces idem-potent iso-mor-phism iso-mor-phisms
multi-plic-a-tion nil-potent poly-nomial priori rami-fication
sin-gu-lar-ities
sub-vari-eties sub-vari-ety trans-form-a-tion trans-form-a-tions
Castel-nuovo
Enri-ques Lo-ba-chev-sky Theo-rem Za-ni-chelli in-vo-lu-tion
Na-ra-sim-han Bohr-Som-mer-feld}

\define\rest#1{_{\textstyle{|}#1}} 

\define\Span#1{\left<#1\right>} 

\define\half{{\textstyle{1\over2}}}


\define\C{\Bbb C} 
\define\R{\Bbb R} 
\define\Z{\Bbb Z} 


\define\sA{{\Cal A}} 
\define\sG{{\Cal G}} 
\define\sM{{\Cal M}} 

\define\al{\alpha}
\define\be{\beta}

\define\om{\omega}

\define\Om{\Omega}









 \document

  \topmatter
  \title The space of hermitian triples and the Ashtekar --- Isham
  quantization
            \endtitle
  \author Nik. Tyurin                     \endauthor

   \address MIIT
   \endaddress
  \email
    tyurin\@tyurin.mccme.ru  ntyurin\@newton.kias.re.kr
   \endemail

\abstract It was pointed out ([5]) that the space of hermitian
triples is an analogy of the hermitian connection space.
Generalizing the Ashtekar --- Isham procedure one can quantize
the space of hermitian triples as well as the original one. Here
we add an example how this similarity can be exploited in a
quantum theory of riemannian geometry.
\endabstract

\endtopmatter

\head \S 0. Introduction   \endhead

Let $X$ be an orientable  compact smooth riemannian 4 -
dimensional manifold. Consider the space of all triples
$$
\sM_X = \{ (g, J, \om) \}
$$
where $g$ is a riemannian metric, $J$ is a compatible almost
complex structure and $\om$ is the corresponding almost kaehler
form. This space naturally splits with respect to a discrete
parametrization such that every connected component can be
described as follows
$$
\sM^K = \{ (g, J, \om) \in \sM_X \quad | \quad K_J = K \in H^2(X,
\Z)\}
$$
where $K_J$ is the canonical class of $J$. All details of the
local geometry can be found in [5], [6].

For the hermitian connection space it was constructed a
quantization programme (see f.e. [2]), dealing with based loops on
3 - dimensional based manifold $M$. One can generalize this
procedure for the space of hermitian triples over the product $X
= M \times S^1$ where for every loop $\al$ on $M$ we have the
corresponding 2 - dimensional torus $\Sigma_{\al}$ on $X$ (or if
one uses graphs for the definition as in [1] then one gets a
Rieman surface). Thus one can arrange a complex bounded function
on the space $\sM^K$ associated with each loop $\al$ and then
repeat the construction of the corresponding $C^*$ - algebra. As
well one can introduce the notion of cylindrical functions and
generate a regular measure on the compactified space
$\overline{\sM^K}$.

Extending this approach to the framework of quantum theory of
riemannian geometry (see [1]) one can exploit "duality" between
two spaces in some special abelian case. In this case we take a 3
- dimensional orientable compact smooth riemannian manifold as
the based manifold and a complex line bundle $L \to M$. The
configuration space is the compactified space of abelian
hermitian connection $\overline{\sA}$ (in  notations of [1]) but
one takes the product
$$
\overline {\sA} \times \overline {\sM^K_M}
$$
as the phase space of the system. Here $\sM^K_M$ is a subspace of
the space $\sM^K$ of all hermitian triples over the product
manifold $X = M \times S^1$ consists of all triples invariant
under all $S^1$ - rotation. The key point is that the product
space
$$
\sA \times \sM^K_M
$$
is a Poisson manifold thus one can say that in some weak sense the
space $\sM^K_M$ is dual to $\sA$. Using this Poisson structure
one could come in the same way as it is proposed for the
cotangent bundle $T^*\sA$ (see section 3A in [1]). Moreover since
the quantizations of both spaces depend on the same loops $\al$
in $M$ one can introduce the notions of cylindrical function
simultaneously on both spaces and then gets well defined Poisson
brackets on the space of cylindrical functions on the product
space.

\head \S 1. The Ashtekar --- Isham quantization
\endhead

In this section we follow [2].

Let $M$ be a real smooth connected orientable 3 - manifold. On
the principal $SU(2)$ - bundle one has the space of hermitian
connections denoted as $\sA$ with natural gauge group $\sG$
action. Fixing a point $x_0 \in M$ consider the space of based
loops denoted as $\Cal L_{x_0}$ which admits a group structure. A
natural equivalence relation on $\Cal L_{x_0}$ is given by the
requirement
$$
H(\al_1, A) = H(\al_2, A)
$$
for every connection $A$ where $H(\al_i, A)$ is the corresponding
holonomy. The collection of the equivalence classes is denoted as
$\Cal H \Cal G$.

Every element $\tilde{\al} \in \Cal H \Cal G$ defines a function
on the quotient space $\sA/ \sG$ by
$$
T_{\tilde{\al} }([A]) = \half Tr H(\al, A)
$$
where $\al$ represents class $\tilde{\al}$ and $A$ represents
class $[A]$. Due to $SU(2)$ trace identities one establishes that
the set of such functions is closed under multiplication thus one
gets a $C^*$ - algebra $\Cal H \Cal A$, consisting of all finite
linear combination of $T_{\tilde{\al}_i}$ with complex
coefficients. This one is a subalgebra in $C^*$ - algebra of all
complex valued bounded continuous functions on $\sA / \sG$ and
one can take the completion $\overline{\Cal H \Cal A}$ of $\Cal H
\Cal A$ under the supremum norm in the ambient $C^*$ - algebra.
This one is called the Ashtekar --- Isham $C^*$ - algebra.

Briefly speaking the Gelfand spectrum of $\overline{\Cal H \Cal
A}$ is a compactification of the quotient space $\sA/\sG$; the
Hilbert space of the quantized theory is represented by
$L^2(\overline{\sA/\sG}, d\mu)$ where $d\mu$ is a regular
diffeomorphism invariant measure on the compactified quotient
space. The construction of such a measure is a crucial step in the
programme. The definition uses the notion of cylindrical functions
--- special functions which form a dense subset in $\Cal H \Cal
A$. Namely for a number of loops $\al_1, ... \al_n$ consider the
following projection
$$
p_{\al_1, ..., \al_n}: \sA/\sG \to (SU(2))^n,
$$
defined by the holonomies around the loops. Then for every
regular function on the target space one takes the corresponding
lifting to the quotient space --- it gives a cylindrical function
on the last one. As well one lifts the product Haar measure to
$\sA/ \sG$; taking the limit one gets a regular $Diff M$ -
invariant measure on the compactified space. Moreover this
measure corresponds to an invariant of framed knots in $M$ (see
[2], [3]).

\head \S 2. The same construction for hermitian triples
\endhead

Let us adopt this construction to the case of hermitian triples.
Instead of sufficiently wide generalization made in [5] now we'd
like to consider only a special case. Namely let $M$ be as above;
consider the product manifold
$$
X = M \times S^1,
$$
and let an orientation on $S^1$ is fixed as well as on $M$. Then
we can consider the space $\sM^K_X$ with a fixed topological datum
$K \in H^2(X, \Z)$. For every loop $\al \in \Cal L_{x_0}$ one has
the  function
$$
P_{\al}: \sM^K_X \to \C
$$
defined as follows. Let $\Sigma_{\al}$ be the torus in $X$
corresponds to $\al$:
$$
\Sigma_{\al} = \al \times S^1.
$$
This torus inherits its own orientation from a parametrization of
$\tilde{\al}$. We take
$$
P_{\al}((g, J, \om)) = e^{\imath(sgn \Sigma_{\al} Vol_g
\Sigma_{\al} + \int_{\Sigma_{\al}} \om)}
$$
where $sgn$ equals to $+1$ if the own orientation of
$\Sigma_{\al}$ is compatible to the fixed orientation on $X$ and
$-1$ otherwise. It's clear that $P_{\al}$ is a bounded continuous
complex function on $\sM^K_X$. Thus again we can make all steps
getting the corresponding $C^*$ - algebra $\overline{\Cal H \Cal
P}$ with Gelfand spectrum denoted as $\overline{\sM^K_X}$. But in
what follows we need just a subspace in $\sM^K_X$ corresponds to
such triples $(g, J, \om) \in \sM^K_X$ which are invariant under
all $S^1$ - rotations defined on $X$ due to its product structure.
This subspace we denote as $\sM_M^K$. It depends only on the
based manifold $M$ and can be described without any references to
the product manifold $X$ (see [7]). Again using the quantization
procedure we get the corresponding compactification
$\overline{\sM^K_M}$.

As above we can define the notion of cylindrical functions on
$\overline{\sM^K_M}$. The definition will be the same: for each
set of loops $\al_1,..., \al_n$ one has a projection
$$
p_{\al_1,...,\al_n} \sM^K_M \to (U(1))^n.
$$
In terms of this projection we again define the cylindrical
functions and the corresponding regular measure. The limit gives
what we need. All details can be found in [7].

\head \S 3. "Duality"
\endhead

In a classical setup the hermitian connection space is the
configuration space. At the same time one takes the cotangent
bundle $T^* \sA$ as the phase space. It is endowed with natural
Poisson brackets but this brackets can not be extended to a
sufficiently wide quantization procedure in the Ashtekar ---
Isham framework (see section 3 A in [1]). We'd like to propose an
appropriate way avoiding the difficulties in the abelian case. We
claim that in this case the space of hermitian triples is almost
dual to the hermitian connection space over a 3 - dimensional
based manifold.

Let $M$ be as above 3 - dimensional based manifold. Consider a
complex line bundle $L$ over $M$ with a fixed hermitian structure.
Let us fix an appropriate element $K \in H^2(X, \Z)$
corresponding to an almost complex structure on $M \times S^1$.
We take the induced component $\sM^K_M$ and consider the direct
product
$$
Y = \sA(L) \times \sM^K_M.
$$

\proclaim{Claim} The space $Y$ is a Poisson manifold.
\endproclaim
(see [8] for the definition and properties of Poisson manifolds
and  [7] for detailed explanations).

The key point is that in each point $(a, (g, J, \om)) \in Y$ the
tangent space is isomorphic to
$$
T_{pt}Y = \Om^1(\imath \R)_M \oplus \Om^2_M \oplus \Om^1_M
$$
and this representation is locally trivial. Thus for any point of
$Y$ the tangent space is represented as
$$
T_{pt} Y = V \oplus V^* \oplus V
$$
(we rescale all vectors from the first summand by $\imath$). For a
pair of smooth functions $F, G \in C^{\infty}(Y \to \R)$ let
$dF_i, dG_i$ be components belong to $i$ - summands in the dual
decomposition
$$
T^*_{pt}Y = V^* \oplus V \oplus V^*.
$$
The Poisson structure is given by formula
$$
\{F, G \}_{\Pi}|_{pt} = dF_1 (dG_2) - d G_1(dF_2) \in \R.
$$
In local coordinate $(p, q, s)$ where $dp$ belongs to the first
summand in the decomposition above $dq$ to the second and $ds$ to
the third the Poisson structure has absolutely classical form
$$
\{F, G\}_{\Pi} = \sum_{i<j}(\frac{\partial F}{\partial p_i}
\frac{\partial G}{\partial q_j} - \frac{\partial F}{\partial q_i}
\frac{\partial G}{\partial p_j}),
$$
while the functions depend only on $s$ form the center of the
corresponding Poisson algebra.

Now we can quantize the picture using the same loops for both
components in $Y$. We introduce the spaces of cylindrical
functions for both $\sA$ and $\sM^K_M$ denoting the first one as
$Cyl_p$ and the second one as $Cyl_q$. Then one has the big space
$Cyl_{tot} = Cyl_p \otimes Cyl_q$ endowed with the product
regular measure. The corresponding limit gives us an universal
space $\Cal C$ on $\overline{Y}$ which is the product of quantized
spaces
$$
\overline{Y} = \overline{\sA} \times \overline {\sM^K_M}
$$
and we have the corresponding $Diff M$ - invariant measure on the
product space. The point is that on this space we have the
extended Poisson brackets
$$
\{ ; \}_{\tilde{\Pi}}: \Cal C \times \Cal C \to \Cal C,
$$
satisfying obviously the distinguishing properties namely
$$
\aligned \{f_1, f_2 \}_{\tilde{\Pi}} \equiv 0, \\
\{F_1, F_2 \}_{\tilde{\Pi}} \equiv 0 \\
\endaligned
$$
where $f_i$ doesn't depend on the second "variable"
$\overline{\sM^K_M}$ being a pure function on $\overline{\sA}$
and $F_i$ depends only on the second "variable".

Consider the  mixed situation when one takes a cylindrical
function $f \in Cyl_p$ and a cylindrical function $F \in Cyl_q$
pairing these functions. Let $f$ is defined by a set of loops
$\al_1, ..., \al_k$ and $F$ uses a set $\be_1, ..., \be_l$. In
the total space $Cyl_{tot}$ we have a component corresponds to the
union set $\al_1, ..., \al_k, \be_1, ... \be_l$. Since the Poisson
structure defined above is a constant Poisson structure then the
Poisson brackets $\{f, F \}_{\Pi}$ should belong to the component
in $Cyl_{tot}$.

Moreover one can consider   the second type functions defined on
$\overline{\sM^K_M}$ as derivations of  the original space
$L^2(\overline{\sA}, d \mu)$. Really if we denote the
corresponding $Diff M$ - invariant measure on the compactified
$\overline{\sM^K_M}$ as $d \mu'$ then the action
$$
N_F(f) = \int_{\overline{\sM^K_M}} \{ f, F \}_{\tilde{\Pi}} d\mu'
$$
is correctly defined and satisfies the Leibnitz rule thus it can
be used in quantization constructions. This action on cylindrical
functions gives  cylindrical functions again but changes a
"grading" because $N_F(f)$ is defined via the extended loop set
$\al_1, ..., \al_k, \be_1, ... \be_l$ rather then the original
function $f$ with the loop set $\al_1, ..., \al_k$. As in [1] the
construction easily extends to the situation with graphs instead
of loop sets used above. The details are contained in [7].

 \head
Conclusion
\endhead

As it was pointed out in [5] one has a twisted version of the
construction above. In the case above our ingredients in the
direct product
$$
Y = \sA \times \sM^K_M
$$
are independent "variables". But one can twist the picture
imposing for example the following condition.  Every hermitian
triple $(g, J, \om)$ from $\sM^K_M$ defines a hermitian structure
on a real rank 2 subbundle of $TM$ (see [5], [7]). So each point
$pt \in \sM^K_M$ gives a complex line bundle with a fixed
hermitian structure. Let us take the corresponding hermitian
connection space $\sA$ as the fiber getting globally topologically
nontrivial bundle
$$
Y \to \sM^K_M
$$
with $\sA$ - fibers. One could try to exploit a Poisson structure
on $Y$ familiar with the original one described above. But there
exists  an alternative way. The point is that $Y$ looks like an
even super symplectic manifold (the definition and properties can
be found in [4]). The based space $\sM^K_M$ has a symplectic form
defined as follows. Let us recall from above that the tangent
space in each point is isomorphic to the direct sum
$$
T_{pt} \sM^K_M \cong \Om^1_M \oplus \Om^2_M.
$$
Thus one has a constant symplectic form defined by
$$
\Om((u_1 \oplus w_1); (u_2 \oplus w_2)) = \int_M (u_1 \wedge w_2 -
u_2 \wedge w_1).
$$
The fiber of $Y$ is an affine space over a real vector space
endowed with  inner product (over point $(g, J, \om)$ the
riemannian metric $g$ should be exploited for the definition).
For the even super  symplectic structure it remains to take an
appropriate affine connection on $Y$. Thus in the construction
one can consider over $Y$ the corresponding super brackets instead
of the classical one.

On the other hand the example has been discussed in the present
article looks like too special. The construction hasn't been
extended to a nonabelian case. The crucial point for a virtual
extension is that the case of hermitian triples is "abelian" in
some sense. One could expect that a nontrivial generalization of
the example lies in the framework of lorentzian geometry  when
one considers triples consist of lorentzian metrics and
"compatible" almost complex structures. Anyway the formula of the
Poisson structure $\{ ; \}_{\Pi}$ on the product space $\sA \times
\sM^K_M$ can be easily generalized to the case when $\sA$ is the
connection space on $SU(2)$- bundle over $M$ but this
generalization doesn't look like a gauge invariant one being the
result of pure formal approach so we finish this example hoping
that more geometric approach will be found in a future.

At the end I'd like to thank Korean Institute for Advanced Study
(Seoul) for hospitality and good working condition.

\enddocument